\documentclass[12pt]{amsart}
\usepackage{latexsym, amsmath, amssymb, amsthm}
\usepackage{epsfig}
\usepackage{amscd}
\usepackage{latexsym}
\usepackage{pstricks,pst-node,cite}

\newtheorem{thm}{Theorem}[section]
\newtheorem{lem}[thm]{Lemma}

\newtheorem{cor}[thm]{Corollary}
\newtheorem{exa}[thm]{Example}

\newtheorem{pro}[thm]{Problem}
\theoremstyle{remark}

\newcommand{\C}[1]{\mathcal #1}

\def \Aut {\mbox{\rm Aut\,}}

\def \Fix{\mbox{\rm Fix\,}}
\def \Top {\mbox{\rm Top\,}}
\def \t {\bar{\tau}}
\newcommand{\comb}[2]{\mbox{$\left(\!\!\begin{array}{c}
            {#1} \\[-0.5ex] {#2} \end{array}\!\!\right)$}}

\newcounter{fignum}
\psset{linecolor=black}\newgray{lightgray}{.9}\newgray{darkgray}{.75}\newgray{pup}{.4}
\psset{linewidth=.6pt,dash=3pt 3pt,doublesep=.05,dotsize=1pt 5}
\SpecialCoor
\newcommand{\middlearrow}{\lput{:U}{\begin{pspicture}[shift=0](0,0)(0,0)
\psline[arrows=->,arrowscale=1.5](2.2pt,0)(2.3pt,0)\end{pspicture}}}

\begin{document}

\title
{Enumerations of finite topologies associated with a finite graph}

\author{Dongseok Kim}
\address{Department of Mathematics \\Kyonggi University
\\ Suwon, 443-760 Korea}
\email{dongseok@kgu.ac.kr}
\thanks{}

\author{Young Soo Kwon}
\address{Department of Mathematics \\Yeungnam University \\Kyongsan, 712-749, Korea}
\email{yskwon@yu.ac.kr}

\author{Jaeun Lee}
\address{Department of Mathematics \\Yeungnam University \\Kyongsan, 712-749, Korea}
\email{julee@yu.ac.kr}
\thanks{}

\maketitle

\begin{abstract}
The number of topologies and non-homeomorphic topologies on a fixed finite set are now known up to $n=18$, $n=16$ but still no complete formula yet~\cite{Sl}.
There are one to one correspondence among topologies, preorder and digraphs.
In this article, we enumerate topologies and non-homeomorphic topologies whose underlying graph is a given finite graph.
\end{abstract}

\bigskip

\section{Introduction}

Topology on a finite set, \emph{finite topology} is often used to
demonstrate interesting phenomena and counterexamples to plausible sounding conjectures.
Finite topology also plays a key role in the theory of image analysis~\cite{Ko1, Ko2}, the structures of molecular~\cite{MS1, MS2}, geometries of finite sets~\cite{RKW}
and digital topology.
One of research on finite topology is to enumerate the number of topologies  and
non-homeomorphic topologies on a fixed finite set, we denoted it by $\t(n)$ and $\hbar(n)$, respectively, where
$n$ is the cardinality of the finite set.
Due to many different works, these numbers are now known up to $n=18$, $n=16$ but still no complete formula yet~\cite{Sl}.

It is not difficult to relate a topology on a finite set with a preorder on the finite set. It can be briefly explained as follows.
Let $X$ be a finite set. A relation $R$ on $X$ is called a \emph{preorder} if it is reflexive and transitive.
For a topology $\C{T}$ on $X$, we define a relation  $R(\C{T})$
on $X$ by a rule that $(x, y)\in R(\C{T})$ if and only if every open set containing $x$ also containing $y$, that is,
$x \in \overline{\{y\}}^\C{T}$, where $\overline{A}^\C{T}$ is the closure of $A$ with respect to the topological $\C{T}$.
Then  $R({\C{T}})$ is a preorder on $X$. We call it \emph{the preorder associated with the topology} $\C{T}$.  Conversely, for  a
preorder $R$ on $X$, let $\C{T}(R)=\{ U \subset X : \forall x \in U, R(x) \subset U \}$, where  $R(x)=\{ y : (x, y)\in R\}$.
Then $\C{T}(R)$ is a topology on $X$ and call it \emph{the topology corresponding to the preorder} $R$.
Notice that $\{R(x) ; x \in X\}$ is a base for the topology $\C{T}(R)$.

Ever since the pioneer work of Evans, Harary and Lynn~\cite{EHL}, counting such topologies can be done
by counting digraphs as follows.
Let $G$ be a finite simple graph with vertex set $V(G)$ and edge set $E(G)$.
We use $|X|$ for the cardinality of a set $X$. For a preorder $R$ on $X$, let $D(R)$ be the direct graph whose vertex set
is $X$ and  arc set is $R\setminus \Delta(X)=\{(x,y) : x\not=y$ and $(x,y)\in R\}$.
For a digraph $D$, \emph{the underlying  graph of} $D$ is the graph whose vertex set is equal to that
of $V(D)$ and edge set is $\{ \{x,y\} : xy$ or $yx$ is an arc of $D \}$.
Notice that there are $3^{|E(G)|}$ directed graphs (or digraphs) whose underlying graph is a given graph $G$.
For a topology $\C{T}$ on $X$, we say the underlying graph of $\C{T}$ is
that of the digraph corresponding to the preorder associated with $\C{T}$. Notice that the digraph $D(R)$
corresponding to a preorder $R$ is transitive, i.e., for any pair of distinct vertices  $a$ and $c$  if $ab$ and $
bc$ are arcs of $D(R)$, then $ac$ is also an arc of $D(R)$. Example~\ref{examp;1v} can help to understand
these relations if the reader is not familiar with graph theory. For terms in graph theory, we refer
to~\cite{GT1}.

There have been a few results on finite topology, preorder and digraphs. Evans $et. al$ found a relation between
labeled topologies on $n$ points and the labeled transitive digraphs with $n$ points~\cite{EHL}. Marijuan found
a few useful properties on digraphs and topologies and relations between finite acyclic transitive digraphs and $\tt{T}_0$ topologies~\cite{Ma}.

In the present article, we mainly focus on graphs and unlabeled(non-homeomorphic) topologies and we
enumerate topologies whose underlying graph is a given finite graph $G$, in other word, we enumerate
the number of topologies and non-homeomorphic topologies with respect to the underlying graphs instead of the number of vertices of graphs.
For a given graph  $G$, let $\Top(G)$ be the set of topologies
having $G$ as its underlying graph. Similarly, let $\t(G)=|\Top(G)|$ and let $\hbar(G)$ be the
number of non-homeomorphic topologies whose underlying graph is $G$. Notice that $\t(G)$ is equal to the number of transitive digraphs whose underlying graph is $G$.

The outline of this article is as follows. First, we will provide
precise definitions and some general formulae in Section~\ref{gen}. In Section~\ref{graph}, we
study how $\t(G), \hbar(G)$ are related with graph operations. In Section~\ref{fixed}, we
find $\t(G), \hbar(G)$ for a few graphs. A brief conclusion and some research problems are presented in Section~\ref{remark}.

\section{General formulae} \label{gen}

For a finite set $X$, let $\Top(X)$ be the set of all topologies on $X$.
 Two topologies $\C{T}_1$ on a set $X$ and $\C{T}_2$ on a set $Y$
are equivalent if the two topological spaces $(X, \C{T}_1)$ and $(X, \C{T}_2)$ are homeomorphic.
For a natural number $n$, let $N_n$ be the set $\{1,2, \ldots,n\}$. For convenience,
let $\t(n)=|\Top(N_n)|$ and let
$\hbar(n)$ be the number of non-homeomorphic topologies on $N_n$.
It is clear that $\t(1)=1$ and $\hbar(1)=1$.

\begin{exa}\label{examp;1v} Let $X=\{a, b\}$ and let  $\C{T}=\{\emptyset, \{a\},  X\}$. Then $\C{T}$ is a topology on $X$.
The preorder $R(\C{T})$ associated with $\C{T}$ is $\{(a,a), (b,b), (b,a)\}$ and the underlying
graph of $\C{T}$ is the complete graph $K_2$ on two vertices $a$ and $b$. Let $R=\Delta(X)=\{(a,a), (b,b))\}$. Then $R$ is a preorder
on $X$. The topology $\C{T}(R)$ associated with $R$ is the discrete topology on $X$ and the underlying
graph of $\C{T}(R)$ is the null graph $\C{N}_2$ on two vertices $a$ and $b$.
\end{exa}

The following example demonstrates the topologies whose underlying graph is the complete graph $K_2$ with $2$ vertices $\{ a, b \}$.

\begin{exa}\label{examp;2v} Let $X=\{a, b\}$. Then  $\C{T}_1=\{\emptyset, X\}$,  $\C{T}_2=\{\emptyset, \{a\}, X\}$, $\C{T}_3=\{\emptyset, \{b\}, X\}$
and $\C{T}_4=\{\emptyset, \{a\}, \{b\}, X\}$ are the list of all topologies on $X$.
So, $\t(2)=4$. Since $(X, \C{T}_2)$ and $(X, \C{T}_3)$ are homeomorphic,
$\hbar(2)=3$. It is clear that  $\Top(K_2)=\{ \C{T}_1, \C{T}_2, \C{T}_3 \}$. This implies that $\t(K_2)=3$ and  $\hbar(K_2)=2$.
\end{exa}

Some properties of finite topological spaces can be
described by graph theoretical terminologies.
The following lemma should be previously known, however, we could not find.
Thus, we provide it in our language.

\begin{lem} \label{lem;propo} Let  $(X, \C{T}_X)$ and   $(Y, \C{T}_Y)$ be finite topological
spaces. Then we have
\begin{enumerate}
\item[{\rm (a)}]  a function  $f: (X, \C{T}_X) \to (Y, \C{T}_Y)$ is continuous if and only if
 $f:(X,R(\C{T}_X)) \to (Y, R(\C{T}_Y))$ preserves the relation, that is,
$f$ is a graph homomorphism between the two digraphs $D(R(\C{T}_X))$ and $D(R(\C{T}_Y))$, 
\item[{\rm (b)}]  the number of components of the topological
space  $(X, \C{T}_X)$ is equal to that of the underlying graph of
$\C{T}_X$. In particular, $(X, \C{T}_X)$ is connected if and only
if the underlying graph of $\C{T}_X$ is connected.
\end{enumerate}
\end{lem}

It comes from Lemma~\ref{lem;propo} (a) that $\hbar(G)$ is equal to
  the number of isomorphism classes of transitive digraphs whose underlying graph is $G$.

Let $\C{T}_1$ and $\C{T}_2$ be two topologies on $N_n$. If they are the same, then their underlying graphs are also the same.
So, two topologies having distinct underlying graphs can not be the same. For a graph $G$, let $\Aut(G)$ be
the group of graph automorphisms of $G$. Then there are $\displaystyle \frac{|V(G)|!}{|\Aut(G)|}$ graphs
 on $V(G)$ that are isomorphic to $G$. The following theorem comes from this observation.

\begin{thm}\label{thm;genfor}
Let $n$ be a natural number. Then we have
$$\t(n)=\sum_G  \frac{n!}{|\Aut(G)|} \t(G)  \mbox{ \rm and } \hbar(n)=\sum_G \hbar(G),$$
where  $G$ runs over all representatives of
isomorphism classes of graphs of $n$ vertices and $\Aut(G)$ is the group of all graph automorphisms of $G$.
\end{thm}

In order to complete the computation $\t(n)$ and $\hbar(n)$, we need to compute $\t(G)$ and $\hbar(G)$ for a given graph $G$.

For a graph  $G$, let $\mathcal{T\!D}(G)$ be the set of transitive digraphs whose underlying graph is $G$. Then
 $\Aut(G)$ acts on the set $\mathcal{T\!D}(G)$ and  $\hbar(G)=|\mathcal{T\!D}(G)/\Aut(G)|$ by Lemma~\ref{lem;propo} (a).
 Now, the following lemma comes from the Burnside lemma.

\begin{lem} \label{lem;burnfor}
Let $G$ be a connected graph and let $\mathcal{T\!D}(G)$ be the set of transitive digraph whose underlying graph is $G$.
Then $$\hbar(G)=\frac{1}{|\Aut(G)|} \sum_{\displaystyle \sigma\in \Aut(G)} |\Fix_\sigma|, $$
where $\Fix_\sigma=\{ D\in \mathcal{T\!D}(G): \sigma(D)=D\}$.
\end{lem}

It is easy to show  that every vertex induced subgraph of a transitive digraph is transitive. From this, we can have the following lemma.

\begin{lem} \label{lem;subgraph}
 For any graph $G$,  $\t(G)=0$ if and only if there exists a vertex induced subgraph $H$ of $G$ such that
  $\t(H)=0$.
\end{lem}

 A graph $G$ is said to be \emph{triangle free} if $G$ dose not contain
any triangles. Let $G$ be a triangle free graph having at most $3$ vertices. Then
every vertex of a transitive digraph having $G$ as its  underlying graph  is a source or a sink.
Since every source is adjacent to a sink and vise versa,
$\t(G)\not=0$ if and only if $G$ is bipartite. To complete the computation $\t(G)$ for a triangle free graph $G$,
we need to compute $\t(G)$ and $\hbar(G)$ for a connected bipartite graph $G$.
If $G$ is a connected bipartite graph having $2$ vertices, then $G$ must be $K_2$. It is observed in Example~\ref{examp;2v}
that $\t(K_2)=3$ and $\hbar(K_2)=2$. Let $G$ be a connected bipartite graph having
at least $3$ vertices and let $X_1$ and $X_2$ be a vertex bipartition of $G$. Let $D$ be a transitive digraph whose underlying graph is $G$.
Then every vertex in $X_1$ is a source and every vertex in $X_2$ is a sink or vice versa.  It implies that
the number of transitive digraphs having $G$ as its underlying graph is $2$. Moreover,  the two digraphs are isomorphic if and only if
$G$ is reflexible. We summarize this discussion as follows.

\begin{thm} \label{thm;bifor}
For a triangle free graph $G$,  $\t(G)\not=0$ if and only if $G$ is bipartite. Moreover,
for a connected bipartite graph $G$ having at least two vertices, we have
 $$\t(G)=\left\{\begin{array}{ll}
 3  & \mbox{ \rm if $G=K_2$}\\
  2  & \mbox{ \rm if $G\not=K_2$}
  \end{array}\right.    \mbox{ \rm and }\,\,
  \hbar(G)=\left\{\begin{array}{ll}
 1  & \mbox{ \rm if $G\not=K_2$ and $G$ is reflexible}\\
  2  & \mbox{ \rm otherwise}
  \end{array}\right.$$
\end{thm}

\section{Topologies and graph operations} \label{graph}

In this section,  we will compute the number $\t(G\circledcirc H)$ and $\hbar(G\circledcirc H)$ when $\t(G)$, $\t(H)$, $\hbar(G)$, and $\hbar(H)$
are known, where $\circledcirc$ is the product, disjoint union, or amalgamation of graphs.

The following lemma gives a computation formula for the graph that can be expressed by a disjoint union of some connected graphs.

\begin{lem} \label{lem;graphfor}
For a natural number $\ell$, let  $G_1,  \ldots, G_{\ell}$  be pairwise nonisomorphic $\ell$ connected graphs and let $n_1, \ldots, n_\ell$ be
$\ell$ natural numbers. Let  $G=n_1G_1 \cup n_2G_2 \cup \cdots \cup n_\ell G_\ell$.
 Then $$\t(G)= \prod_{i=1}^\ell \t(G_i)^{n_i}\quad \mbox{\rm  and} \quad
 \hbar(G)= \prod_{i=1}^\ell \hbar(n_iG_i)=\prod_{i=1}^\ell \comb{\hbar(G_i)+n_i-1}{n_i},$$
 where  $mH$ stands for the disjoint union of $m$ copies of $H$.
\begin{proof}
Let  $H_1$ and $H_2$ be two graphs having disjoint vertex sets. Then $\t(H_1 \cup H_2)=\t(H_1)\t(H_2)$.
It implies that $\t(G)= \prod_{i=1}^\ell \t(G_i)^{n_i}$. To prove the second statement, we firstly show that
$\hbar(nH)=\comb{\hbar(H)+n-1}{n}$.
Let us identify a topology $\C{T}\in \Top(nH)$ with a finite sequence $(\C{T}_1, \C{T}_2, \ldots, \C{T}_n)$ of length $n$ with
$\C{T}_i\in \Top(H)$ for each $i=1,2,\ldots, n$. Then two topologies $(\C{T}_1, \C{T}_2, \ldots, \C{T}_n)$ and
$(\C{T}_1', \C{T}_2', \ldots, \C{T}_n')$ are equivalent if and only if there exists a bijection $\sigma:N_n \to N_n$ such that
$\C{T}_i$ and $\C{T}_{\sigma(i)}'$ are equivalent for each $i=1,2,\ldots, n$. It implies that the number  $\hbar(nH)$ is equal to
the number of selections with repetitions of $n$ objects chosen from $\hbar(H)$ types of objects, i.e.,
$\hbar(nH)=\comb{\hbar(H)+n-1}{n}$. For given two nonisomorphic graphs
 $H_1$ and $H_2$, it is not hard to show that  $\hbar(H_1 \cup H_2)=\hbar(H_1)\hbar(H_2)$. It completes the proof.
\end{proof}
\end{lem}

For two graphs $G$ and $H$, the Cartesian product $G \Box H$ is a
graph such that $V(G \Box H)=V(G) \times V(H)$ and two vertices
$(u_1, v_1)$ and $(u_2, v_2)$  are adjacent  if and only if
($u_1=u_2$ and $v_1v_2 \in E(H)$) or ($u_1u_2 \in V(G)$ and
$v_1=v_2$). We aim to compute $\t(G \Box H)$ and $\hbar(G \Box
H)$.

\begin{lem} \label{K2-cart-oddcycle}
For any odd number $n\ge 3$,  $\t(K_2 \Box C_n)=0$.
\begin{proof}
 If $n$ is greater than or equal to $5$, then it is not hard to show that
 there is no transitive digraph whose underlying graph is $C_n$.
Since $K_2 \Box C_n$ has an induced subgraph isomorphic to $C_n$, $\t(K_2 \Box C_n)=0$ by Lemma~\ref{lem;subgraph}.

Let $n=3$ and let $D$ be a digraph whose underlying graph is $K_2
\Box C_3$. Let $V(K_2) = \{ u_1 , u_2 \}$ and $V(C_3 ) = \{ v_1 ,
v_2 , v_3 \}$. One can check that there are $u_i \in V(K_2)$ and
$v_j, v_k \in V(C_3)$ such that both $((u_i, v_j), (u_{3-i},
v_j))$ and $((u_{3-i}, v_j), (u_{3-i}, v_k))$ are directed edges
in $D$. Since $((u_i, v_j)$ and $(u_{3-i}, v_k))$ are not adjacent in $K_2
\Box C_3$, $D$ is not transitive.  Hence $\t(K_2 \Box
C_3)=0$.
\end{proof}
\end{lem}

Note that for two graphs $G$ and $H$, $G \Box H$ is bipartite if
and only if both $G$ and $H$ are bipartite. Furthermore, for
bipartite graph $G \Box H$, one can check that  $G \Box H$ is
reflexible if and only if either $G$ or $H$ is reflexible.

\begin{thm} \label{cartesian-product}
 For any nontrivial graphs $G$ and $H$,
 $$\t(G \Box H)=\left\{\begin{array}{ll}
 0  & \mbox{ \rm if $G$ or $H$ is not bipartite,}\\
  2  & \mbox{ \rm otherwise}
  \end{array}\right.$$ and
 $$\hbar(C_n)=\left\{\begin{array}{ll}
 0  & \mbox{ \rm if $G$ or $H$ is not bipartite,}\\
  1  & \mbox{ \rm if (both $G$ and $H$ are bipartite) and (either $G$ or $H$ is reflexible),}\\
  2  & \mbox{ \rm otherwise.}
  \end{array}\right.$$
\begin{proof}
Assume that $G$ or $H$ is not bipartite. Then $G \Box H$ is not
bipartite and $G \Box H$ contains an induced subgraph isomorphic
to $K_2 \Box C_n$ for some odd $n$. By Lemmas~\ref{lem;subgraph}
and Lemma~\ref{K2-cart-oddcycle}, we have $\t(G \Box H)=0$.

Suppose $G$ and $H$ are bipartite. Now $G \Box H$ is bipartite and
hence $\t(G \Box H)=2$ by Theorem~\ref{thm;bifor}. Furthermore,
$\hbar(G \Box H)$ is $1$ if both $G$ and $H$ are reflexible; 2
otherwise.
\end{proof}
\end{thm}

Let $v$ be a cut vertex of a graph $G$. Then $v$ is a sink or a source in every transitive digraph having $G$ as its underlying graph.
Let $G$ and $H$ be two graphs. For two vertices $u\in V(G)$ and  $v \in V(H)$, a graph $G*_{u=v}H$ is obtained from $G$ and $H$ by identifying the vertices $u$ and $v$. We call it \emph{the amalgamation of $G$ and $H$ along the vertices $u$ and $v$}.
We note that $u$ (or $v$) is a cut vertex of $G*_{u=v}H$. For a graph $G$ and a vertex $u\in V(G)$,
let $\t_{si}(G, u)$ ($\t_{so}(G,u)$, respectively) be the number of
topologies having $G$ as its underlying graph and  $v$ as a sink(source, respectively) in the digraphs corresponding them.
Similarly, we define $\hbar_{si}(G,u)$ and $\hbar_{so}(G,u)$. For any transitive digraph $D$,
if we change direction of every directed edge in $D$, then the resulting digraph is also a transitive digraph. Moreover the correspondence is 
bijective and hence $\t_{si}(G, u) = \t_{so}(G,u)$ and $\hbar_{si}(G,u)= \hbar_{so}(G,u)$.
Note that for any preorder $R$, the digraph obtained by changing direction of every edge in $D(R)$ corresponds to
the topology composed of all closed sets in $\C{T}(R)$.

Now the following lemma comes from the fact that $v$ is a cut vertex of $G*_{u=v}H$.
\begin{lem} \label{cut-vertex} Let $G$ and $H$ be two graphs.
\begin{enumerate}
\item[{\rm (a)}] If $v$ is a cut vertex of $G$, then $v$ is sink or source in every transitive digraph having $G$ as its underlying graph and hence $\t(G)=2\t_{si}(G,v)$. Moreover, if $v$ is the unique cut vertex, then  $\hbar(G)=2\hbar_{si}(G,v).$
\item[{\rm (b)}]If  $u\in V(G)$ and  $v \in V(H)$, then
 $\t(G*_{u=v}H)=2\t_{si}(G, u)\t_{si}(H, v).$
 Moreover, both either $G$ and $H$ have no cut vertices, then
 $$\hbar(G*_{u=v}H)=\left\{\begin{array}{ll}
 (\hbar_{si}(G, u)+1)\, \hbar_{si}(G, u)&
  \mbox{\rm if $u\simeq_fv$,}\\[2ex]
 2\hbar_{si}(G, u)\,\hbar_{si}(H, v) &
  \mbox{\rm otherwise},
 \end{array}  \right.$$
 where $u\simeq_fv$ means that there exists a graph isomorphism $f:G\to H$ such that
 $f(u)=v$.
 \end{enumerate}
\end{lem}

From Theorem~\ref{thm;genfor} and Lemma~\ref{lem;graphfor}, we can see that the computation of  $\t(n)$ and $\hbar(n)$ can be completed
 if we can compute $\t(G)$ and $\hbar(G)$ for any connected graph with $n$ vertices.
In the next section, we compute $\t(G)$ and $\hbar(G)$
 for some special classes of connected graphs.

\section{Topologies having a fixed underlying graph} \label{fixed}

In this section, we compute $\t(G)$ and $\hbar(G)$ when $G$ is a cycle, a wheel or a complete graph.

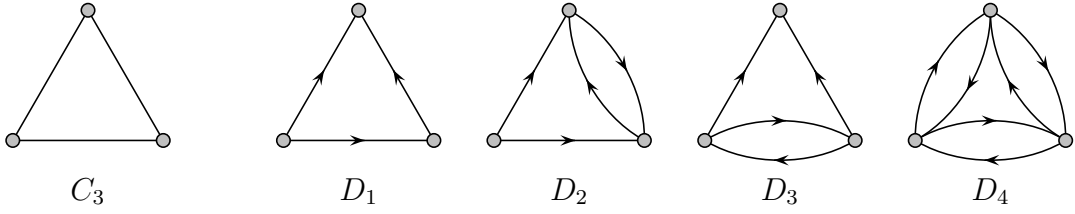
\begin{figure}
$$
\begin{pspicture}[shift=-2](-.3,-1)(2.3,2) \psline(0,0)(1,1.732)
\psline(0,0)(2,0) \psline(1,1.732)(2,0)
\pscircle[fillstyle=solid,fillcolor=darkgray](0,0){.1}
\pscircle[fillstyle=solid,fillcolor=darkgray](1,1.732){.1}
\pscircle[fillstyle=solid,fillcolor=darkgray](2,0){.1}
\rput[t](1,-.5){$C_3$}
\end{pspicture}
\hskip 1cm
\begin{pspicture}[shift=-1.9](0,-1)(0,2)
\end{pspicture}\begin{pspicture}[shift=-1.9](0,2)(0,2)
\begin{pspicture}[shift=-1.9](-.3,-1)(2.3,2)
\rput(0,0){\rnode{a1}{$$}}
\rput(1,1.732){\rnode{a2}{$$}}
\rput(2,0){\rnode{a3}{$$}}
\ncline{a1}{a3}\middlearrow
\ncline{a1}{a2}\middlearrow \ncline{a3}{a2}\middlearrow
\pscircle[fillstyle=solid,fillcolor=darkgray](0,0){.1}
\pscircle[fillstyle=solid,fillcolor=darkgray](1,1.732){.1}
\pscircle[fillstyle=solid,fillcolor=darkgray](2,0){.1}
\rput[t](1,-.5){$D_1$}
\end{pspicture}\end{pspicture}
\hskip .2cm \begin{pspicture}[shift=-1.9](0,2)(0,2)
\begin{pspicture}[shift=-1.9](-.3,-1)(2.3,2)
\rput(0,0){\rnode{a1}{$$}}
\rput(1,1.732){\rnode{a2}{$$}}
\rput(2,0){\rnode{a3}{$$}}
\ncline{a1}{a3}\middlearrow
\ncline{a1}{a2}\middlearrow
\nccurve[angleA=-30,angleB=90,nodesep=1pt]{a2}{a3}\middlearrow
\nccurve[angleA=150,angleB=-90,nodesep=1pt]{a3}{a2}\middlearrow
\pscircle[fillstyle=solid,fillcolor=darkgray](0,0){.1}
\pscircle[fillstyle=solid,fillcolor=darkgray](1,1.732){.1}
\pscircle[fillstyle=solid,fillcolor=darkgray](2,0){.1}
\rput[t](1,-.5){$D_2$}
\end{pspicture}\end{pspicture}
\hskip .2cm 
\begin{pspicture}[shift=-1.9](0,2)(0,2)
\begin{pspicture}[shift=-1.9](-.3,-1)(2.3,2)
\rput(0,0){\rnode{a1}{$$}}
\rput(1,1.732){\rnode{a2}{$$}}
\rput(2,0){\rnode{a3}{$$}}
\nccurve[angleA=30,angleB=150,nodesep=1pt]{a1}{a3}\middlearrow
\nccurve[angleA=-150,angleB=-30,nodesep=1pt]{a3}{a1}\middlearrow
\ncline{a1}{a2}\middlearrow \ncline{a3}{a2}\middlearrow
\pscircle[fillstyle=solid,fillcolor=darkgray](0,0){.1}
\pscircle[fillstyle=solid,fillcolor=darkgray](1,1.732){.1}
\pscircle[fillstyle=solid,fillcolor=darkgray](2,0){.1}
\rput[t](1,-.5){$D_3$}
\end{pspicture}\end{pspicture}
\hskip .2cm\begin{pspicture}[shift=-1.9](0,2)(0,2)
\begin{pspicture}[shift=-1.9](-.3,-1)(2.3,2)
\rput(0,0){\rnode{a1}{$$}}
\rput(1,1.732){\rnode{a2}{$$}}
\rput(2,0){\rnode{a3}{$$}}
\nccurve[angleA=-30,angleB=90,nodesep=1pt]{a2}{a3}\middlearrow
\nccurve[angleA=150,angleB=-90,nodesep=1pt]{a3}{a2}\middlearrow
\nccurve[angleA=30,angleB=150,nodesep=1pt]{a1}{a3}\middlearrow
\nccurve[angleA=-150,angleB=-30,nodesep=1pt]{a3}{a1}\middlearrow
\nccurve[angleA=90,angleB=-150,nodesep=1pt]{a1}{a2}\middlearrow
\nccurve[angleA=-90,angleB=30,nodesep=1pt]{a2}{a1}\middlearrow
\pscircle[fillstyle=solid,fillcolor=darkgray](0,0){.1}
\pscircle[fillstyle=solid,fillcolor=darkgray](1,1.732){.1}
\pscircle[fillstyle=solid,fillcolor=darkgray](2,0){.1}
\rput[t](1,-.5){$D_4$}
\end{pspicture}\end{pspicture}
$$
\caption{$C_3$ and four nonisomorphic transitive digraphs whose underlying graph is $C_3$.} \label{fig1}
\end{figure}

First, we aim to compute $\t(C_n)$ and $\hbar(C_n)$ for a natural number $n \ge 3$.
Notice that $C_3$ is the complete graph $K_3$ on three vertices. There are
$4$ nonisomorphic transitive digraphs whose underlying graph is $K_3$  as illustrated in Figure~\ref{fig1} and the following is the list of all representatives of them;
 $D_1=\{ 1 | 2 | 3\}$,
$D_2=\{ 1 | 2, 3 \}$, $D_3=\{ 1, 2 | 3 \}$,
$D_4=\{ 1, 2, 3 \}$, where $D_3=\{ 1, 2 | 3 \}$ stands for the digraph with vertex
set $\{1,2,3\}$ and arc
set  $\{12, 21, 13, 23\}$. Hence $\hbar(C_3)=4$. For
convenience, let $\alpha_i$ be the number of digraphs that are
isomorphic to $D_i$, Then $\alpha_1=4!=6$,
$\alpha_2=\alpha_3=\dfrac{3!}{2!}=3$,
$\alpha_4=\dfrac{3!}{3!}=1$. Hence,
$\t(C_3)=\sum_{i=1}^8\alpha_i=6+6+1=13$. Since
 there is no transitive digraphs whose underlying graph is $C_n$ when $n$ is odd greater than $3$.

Now, the following corollary comes from Lemma~\ref{lem;subgraph}, Theorem~\ref{thm;bifor} and the fact that every  cycle of even length is reflexible.

\begin{cor}\label{cor;cycle}
 For a natural number $n\ge 3$,
$$\t(C_n)=\left\{\begin{array}{ll}
13  & \mbox{ \rm if $n=3$},\\
 0  & \mbox{ \rm if $n$ is odd and $n \ge 5$},\\
  2  & \mbox{ \rm otherwise},
  \end{array}\right.$$
   and
   $$
  \hbar(C_n)=\left\{\begin{array}{ll}
4  & \mbox{ \rm if $n=3$},\\
 0  & \mbox{ \rm if $n$ is odd and $n \ge 5$},\\
  1  & \mbox{ \rm otherwise}.
  \end{array}\right.$$
\end{cor}

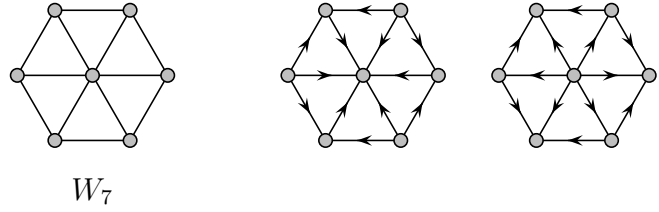
\begin{figure}
$$
\begin{pspicture}[shift=-2](-.3,-.8)(2.3,2) \psline(0,.866)(2,.866)
\psline(.5,1.732)(1.5,1.732) \psline(.5,0)(1.5,0)
\psline(0,.866)(.5,1.732)
\psline(.5,0)(1.5,1.732) \psline(1.5,0)(2,.866)
\psline(0,.866)(.5,0)
\psline(.5,1.732)(1.5,0) \psline(1.5,1.732)(2,.866)
\pscircle[fillstyle=solid,fillcolor=darkgray](0,.866){.1}
\pscircle[fillstyle=solid,fillcolor=darkgray](.5,1.732){.1}
\pscircle[fillstyle=solid,fillcolor=darkgray](.5,0){.1}
\pscircle[fillstyle=solid,fillcolor=darkgray](1,.866){.1}
\pscircle[fillstyle=solid,fillcolor=darkgray](1.5,1.732){.1}
\pscircle[fillstyle=solid,fillcolor=darkgray](1.5,0){.1}
\pscircle[fillstyle=solid,fillcolor=darkgray](2,.866){.1}
\rput[t](1,-.5){$W_7$}
\end{pspicture}
\hskip 1cm
\begin{pspicture}[shift=-1.9](0,-.8)(0,2)
\end{pspicture}\begin{pspicture}[shift=-1.9](0,2)(0,2)
\begin{pspicture}[shift=-1.9](-.3,-.8)(2.3,2)
\rput(1,.866){\rnode{a0}{$$}}
\rput(0,.866){\rnode{a1}{$$}}
\rput(.5,0){\rnode{a2}{$$}}
\rput(1.5,0){\rnode{a3}{$$}}
\rput(2,.866){\rnode{a4}{$$}}
\rput(1.5,1.732){\rnode{a5}{$$}}
\rput(.5,1.732){\rnode{a6}{$$}}
\ncline{a1}{a0}\middlearrow
\ncline{a2}{a0}\middlearrow
\ncline{a3}{a0}\middlearrow
\ncline{a4}{a0}\middlearrow
\ncline{a5}{a0}\middlearrow
\ncline{a6}{a0}\middlearrow
\ncline{a1}{a2}\middlearrow \ncline{a1}{a6}\middlearrow
\ncline{a3}{a2}\middlearrow \ncline{a5}{a6}\middlearrow
\ncline{a3}{a4}\middlearrow \ncline{a5}{a4}\middlearrow
\pscircle[fillstyle=solid,fillcolor=darkgray](0,.866){.1}
\pscircle[fillstyle=solid,fillcolor=darkgray](.5,1.732){.1}
\pscircle[fillstyle=solid,fillcolor=darkgray](.5,0){.1}
\pscircle[fillstyle=solid,fillcolor=darkgray](1,.866){.1}
\pscircle[fillstyle=solid,fillcolor=darkgray](1.5,1.732){.1}
\pscircle[fillstyle=solid,fillcolor=darkgray](1.5,0){.1}
\pscircle[fillstyle=solid,fillcolor=darkgray](2,.866){.1}
\end{pspicture}\end{pspicture}
\hskip .2cm \begin{pspicture}[shift=-1.9](0,2)(0,2)
\begin{pspicture}[shift=-1.9](-.3,-.8)(2.3,2)
\rput(1,.866){\rnode{a0}{$$}}
\rput(0,.866){\rnode{a1}{$$}}
\rput(.5,0){\rnode{a2}{$$}}
\rput(1.5,0){\rnode{a3}{$$}}
\rput(2,.866){\rnode{a4}{$$}}
\rput(1.5,1.732){\rnode{a5}{$$}}
\rput(.5,1.732){\rnode{a6}{$$}}
\ncline{a0}{a1}\middlearrow
\ncline{a0}{a2}\middlearrow
\ncline{a0}{a3}\middlearrow
\ncline{a0}{a4}\middlearrow
\ncline{a0}{a5}\middlearrow
\ncline{a0}{a6}\middlearrow
\ncline{a1}{a2}\middlearrow \ncline{a1}{a6}\middlearrow
\ncline{a3}{a2}\middlearrow \ncline{a5}{a6}\middlearrow
\ncline{a3}{a4}\middlearrow \ncline{a5}{a4}\middlearrow
\pscircle[fillstyle=solid,fillcolor=darkgray](0,.866){.1}
\pscircle[fillstyle=solid,fillcolor=darkgray](.5,1.732){.1}
\pscircle[fillstyle=solid,fillcolor=darkgray](.5,0){.1}
\pscircle[fillstyle=solid,fillcolor=darkgray](1,.866){.1}
\pscircle[fillstyle=solid,fillcolor=darkgray](1.5,1.732){.1}
\pscircle[fillstyle=solid,fillcolor=darkgray](1.5,0){.1}
\pscircle[fillstyle=solid,fillcolor=darkgray](2,.866){.1}
\end{pspicture}\end{pspicture}
$$
\caption{$W_7$ and two nonisomorphic transitive digraphs whose underlying graph is $W_7$.} \label{fig2}
\end{figure}

Next, we will compute $\t(G)$ when $G$ is the wheel graph.
For a  natural number $n \ge 4$, the wheel $W_{n}$  is  a  graph
with  $n$ vertices which contains a cycle $C_{n-1}$ as an induced
subgraph, and every  vertex in the cycle is adjacent to one other
 vertex. Note that $W_4$ is the complete graph $K_4$
on $4$ vertices.
A wheel graph $W_7$ and $4$ nonisomorphic transitive digraphs whose underlying graph is $W_7$ is given in Figure~\ref{fig2}.
There are $8$ nonisomorphic transitive digraphs whose underlying graph is
$K_4$ and the following is the list of all representatives of them;
 $D_1=\{ 1| 2 | 3 | 4 \}$,
$D_2=\{ 1 | 2 | 3, 4 \}$, $D_3=\{ 1 | 2, 3 | 4\}$, $D_4=\{ 1, 2 | 3 | 4 \}$, $D_5=\{ 1 | 2, 3, 4
\}$, $D_6=\{ 1, 2, 3 | 4 \}$, $D_7=\{ 1, 2 | 3, 4 \}$,
and $D_8=\{ 1, 2, 3, 4 \}$.
For
convenience, let $\beta_i$ be the number of digraphs that are
isomorphic to $D_i$, Then $\beta_1=24$,
$\beta_2=\beta_3=\beta_4=12$,
$\beta_5=\beta_6=4$, $\beta_7=6$, and $\beta_8=1$. Hence
$\t(K_4)=\sum_{i=1}^8\beta_i=24+36+8+6+1=75$.
 Now, the following comes from Lemma~\ref{lem;subgraph} and a simple computation.

\begin{thm}\label{thm;wheel}
 For a natural number $n\ge 4$, let $W_n$ be the wheel graph. Then we have
$$
\t(W_n)=\left\{\begin{array}{ll}
4 & \mbox{\rm if $n$ is odd and $n\ge 7$},\\
0 & \mbox{\rm if $n$ is even and $n\ge 6$},\\
75 & \mbox{\rm if $n=4$},\\
8 & \mbox{\rm if $n=5$},\end{array}\right.$$ and
$$\hbar(W_n)=\left\{\begin{array}{ll}
2 & \mbox{\rm $n$ is odd and $n\ge 7$},\\
0 & \mbox{\rm $n$ is even and $n\ge 6$},\\
8 & \mbox{\rm if $n=4$},\\
4 & \mbox{\rm if $n=5$}.\end{array}\right.
$$
\end{thm}

Finally, we will compute $\t(K_n)$ and $\hbar(K_n)$ for the complete graph $K_n$ on $n$ vertices.

\medskip

\begin{thm} \label{thm;k_nfor}
 For a natural number $n$,  we have
 $$\t(K_n)=\sum_{k=1}^n Surj(n,k)=\sum_{k=1}^n S(n,k) k! \quad\mbox{ \rm and }\quad \hbar(K_n)=2^{n-1},$$
 where $Surj(n,k)$ is the number of surjections from $N_n$ to $N_k$ and $S(n,k)$ is the number of ways
 of partitions of $N_n$ into exactly $k$ nonempty parts which is known as a \emph{Stirling number
 the second kind}.
\begin{proof}
Let $R$ be a preorder on $N_n$ whose underlying graph is the
complete graph $K_n$. Then $(i,j)\in R$ or $(j,i)\in R$ for any
two distinct elements  $i$ and $j$ in $N_n$. We define another
relation $E(R)$ on $N_n$ by $(x,y)$ in  $E(R)$  if and only if
both $(x,y)$ and $(y,x)$ are in $R$. Then $E(R)$  is an
equivalence relation on $N_n$. Let  $f_R$ be a relation on
$N_n/E(R)$ defined by $([x],[y]) \in \widetilde{R}$ if and only if
$(x,y) \in R$ for some $x \in [x]$ and $y \in [y]$. By
transitivity, if $(x,y) \in R$ for some $x \in [x]$ and $y \in
[y]$ then $(x',y') \in R$ for all $x' \in [x]$ and $y' \in [y]$.
Now $\widetilde{R}$ is a total order. For convenience, let
$N_n/E(R)=\{ [i_1], [i_2], \ldots, [i_k]\}$ and $([i_s], [i_t])\in
\widetilde{R}$ if and only if $s \le t$. We define $f_R: N_n\to
N_k$ by $f_R(i)=s$ if $(i, i_s)\in E(R)$. Then $f_R$ is a
surjection. Conversely, for a given surjection $f: N_n \to N_k$ we
define a relation $R_f$ on $N_n$ by $(i, j)\in R_f$ if and only if
$f(i) \le f(j)$. Then $R_f$ is a preorder on $N_n$ whose
underlying graph is $K_n$. Since $R_{f_R}=R$ and $f_{R_f}=f$, the
correspondence is one-to-one and hence $\t(K_n)=\sum_{k=1}^n
Surj(n,k)$. Since $Surj(n,k)=S(n,k) k!$, $\t(K_n)=\sum_{k=1}^n
Surj(n,k)=\sum_{k=1}^nS(n,k) k!$.

For a proof of the second equation, let $f:N_n\to N_h$
and $g:N_n\to N_k$ be two surjections. Then $\C{T}(R_f)$ and
$\C{T}(R_r)$ are equivalent if and only if $h=k$ and
$|f^{-1}(i)|=|g^{-1}(i)|$ for all $i=1,2, \ldots, h=k$ (by
Lemma~\ref{lem;propo}). For  a $k$-tuple  $(n_1, n_2, \ldots,
n_k)$ of natural numbers such that $n=n_1+n_2+ \cdots +n_k$, we
define a surjection $\varphi:N_n\to N_k$ such that
$\varphi^{-1}(1)=\{1, 2, \cdots, n_1\}$ and
$$\varphi^{-1}(i)=\left\{\left(\sum_{t=1}^{i-1}n_t\right) +1, \left(\sum_{t=1}^{i-1}n_t\right) +2, \ldots, \left(\sum_{t=1}^{i-1}n_t\right)+n_i-1, \left(\sum_{t=1}^{i}n_t\right)\right\}$$
for each $i=2,\ldots,k$.
 Then the topology corresponding to $\varphi$ is a representative of the equivalence class of topologies corresponding to all surjection $f$
 satisfying $|f^{-1}(i)|=n_i$ for all $i=1,2,\ldots,k$. It is clear that two  different $k$-tuples represent two different equivalence topologies.
 So, the number of equivalence classes of topologies corresponding to the set of all surjections from $N_n \to N_k$ is equal to the number of ways
 to choose $k-1$ positions among $n-1$ positions between the $n$ numbers $1,2, \ldots, n$. Hence we have
$\displaystyle \hbar(K_n)=\sum_{k=1}^{n} \comb{n-1}{k-1}=\sum_{k=0}^{n-1} \comb{n-1}{k}=2^{n-1}$.
\end{proof}
\end{thm}

\section{Further remarks} \label{remark}

We already know that $\t(1)=1$, $\hbar(1)=1$, $\t(2)=4$ and $\hbar(2)=3$. In order to compute
$\t(3)$ and $\hbar(3)$, we list all representatives of isomorphism classes of graphs on $3$ vertices as follows;
 the null graph $\C{N}_3$, $H$,
the path $P_2$ of length 2, and
the complete graph  $K_3$, where $H$ is the disjoint union of $K_2$ and the null graph $\C{N}_1$.
It is not hard to show that $\t(\C{N}_3)=1$, $\hbar(\C{N}_3)=1$, $\t(H)=3$, $\hbar(H)=2$, $\t(P_2)=2$, $\hbar(P_2)=2$, $\t(K_3)=13$, and $\hbar(K_3)=4$.
Now, It comes from Theorem~\ref{thm;genfor} that $\t(3)=1+3 \times 3 + 2 \times 3 + 13=29$
and $\hbar(3)=1+ 2 + 2 + 4 = 9$. Similarly, we can see that $\t(4)=355$ and $\hbar(4)=33$.

At last, we conclude the article by proposing several further research problems as follows.

\begin{pro}
Let $G, H$ be graphs and $K_n$ be a complete graph with $n$ vertices.

\begin{enumerate}
\item[{\rm (a)}] Compute $\t(G)$ and $\hbar(G)$ for arbitrary graph $G$.
\item[{\rm (b)}] Compute $\t(G \circledcirc H)$ and $\hbar(G\circledcirc H)$ for any two distinct graphs
and any binary operation $\circledcirc$ on the set of graphs.
\item[{\rm (c)}] Find an explicit numerical formula to compute $\t(K_n)$.
\item[{\rm (d)}] Find an explicit numerical formula to compute $\t(n)$ and $\hbar(n)$.
\end{enumerate}
\end{pro}

\section*{Acknowledgments}
The \TeX\, macro package
PSTricks~\cite{PSTricks} was essential for typesetting the equations
and figures.

\end{document}